\documentclass[a4paper,12pt]{article}

\usepackage{amsthm}
\usepackage{amsmath,amssymb,latexsym,amsfonts}

\theoremstyle{plain}
\newtheorem{Thm}{Theorem}[section]
\newtheorem{Lem}[Thm]{Lemma}
\newtheorem{Prop}[Thm]{Proposition}

\theoremstyle{definition}

\newtheorem{Rem}[Thm]{Remark}

\newcommand{\Proof}[2][Proof]{\begin{proof}[{#1}] #2 \end{proof}}

\renewcommand{\Re}{\mathop{\rm Re}}
\renewcommand{\Im}{\mathop{\rm Im}}
\renewcommand{\d}{{\rm d}} 
\newcommand{\sgn}{\mathop{\rm sgn}}

\newcommand{\e}{{\rm e}} 

\newcommand{\law}{\stackrel{{\rm law}}{=}}

\newcommand{\tend}[2]{\mathrel{\mathop{\longrightarrow}\limits^{#1}_{#2}}}

\renewcommand{\bar}{\overline}
\newcommand{\un}[1]{\underline{#1}}

\numberwithin{equation}{section}

\makeatletter
\renewcommand\section{\@startsection {section}{1}{\z@}%
                                   {-3.5ex \@plus -1ex \@minus -.2ex}%
                                   {2.3ex \@plus.2ex}%
                                   {\normalfont\large\bf}}
\makeatother

\makeatletter
\renewcommand\subsection{\@startsection {subsection}{1}{\z@}%
                                   {-3.5ex \@plus -1ex \@minus -.2ex}%
                                   {2.3ex \@plus.2ex}%
                                   {\normalfont\normalsize\bf}}
\makeatother

\newcommand{\absol}[1]{\left| #1 \right|} 
\newcommand{\rbra}[1]{\!\left( #1 \right)} 
\newcommand{\cbra}[1]{\!\left\{ #1 \right\}} 
\newcommand{\sbra}[1]{\!\left[ #1 \right]} 

\newcommand{\bC}{\ensuremath{\mathbb{C}}}
\newcommand{\bD}{\ensuremath{\mathbb{D}}}
\newcommand{\bE}{\ensuremath{\mathbb{E}}}
\newcommand{\bP}{\ensuremath{\mathbb{P}}}
\newcommand{\bR}{\ensuremath{\mathbb{R}}}
\newcommand{\cB}{\ensuremath{\mathcal{B}}}
\newcommand{\cF}{\ensuremath{\mathcal{F}}}
\newcommand{\vn}{\ensuremath{\mbox{{\boldmath $n$}}}}

\begin{document}

\begin{center}
{\large \bf 
On harmonic function for the killed process upon hitting zero 
of asymmetric L\'evy processes 
} 
\end{center}
\begin{center}
Kouji Yano\footnote{Graduate School of Science, Kyoto University. 
Email: {\tt kyano@math.kyoto-u.ac.jp}} 
\end{center}
\bigskip

\begin{abstract}
For a certain class of asymmetric L\'evy processes 
where the origin is regular for itself, 
the renormalized zero resolvent is proved to be harmonic 
for the killed process upon hitting zero. 
\end{abstract}


\section{Introduction}

Let $ \{ X(t),\bP_x \} $ be the canonical representation 
of a reflecting Brownian motion on $ [0,\infty ) $ 
and let $ \{ X(t),\bP^0_x \} $ denote its killed process upon hitting zero. 
Let $ \vn $ denote the characteristic measure of excursions away from zero. 
It is then well-known that 
\begin{align}
\bE^0_x[X(t)] = x 
\quad \text{for all $ t>0 $ and $ x \ge 0 $} 
\label{}
\end{align}
and 
\begin{align}
\text{$ \vn[X(t)] $ is finite constant in $ t>0 $}, 
\label{}
\end{align}
where we write $ \vn[F] $ for $ \int F \d \vn $. 
In this case, we say that the function $ h(x) \equiv x $ is 
{\em harmonic for the killed process}. 
Moreover, it is well-known that 
the harmonic transform $ \{ X(t),\bP^h_x \} $ 
of the killed process with respect to $ h $ 
is the 3-dimensional Bessel process, 
and that the excursion measure $ \vn $ enjoys the disintegration formula 
\begin{align}
\vn(\cdot) = \int_0^{\infty } \bP^h_0(\cdot | X(t)=0) \rho(t) \d t , 
\label{eq: BM disinteg}
\end{align}
where $ \rho(t) $ is a certain non-negative function. 

The author \cite{MR2603019} has obtained 
disintegration formula similar to \eqref{eq: BM disinteg} 
for one-dimensional symmetric L\'evy processes, 
where the renormalized zero resolvent 
was proved to be harmonic for the killed process. 
In this paper, we study the asymmetric case. 

This paper is organized as follows. 
In Section \ref{sec: G}, we recall, without proofs, several theorems 
of the excursion theory for L\'evy processes. 
In Section \ref{sec: R}, we study the renormalized zero resolvent. 
In Section \ref{sec: S}, we recall the earlier results in the symmetric case. 
In Section \ref{sec: SS}, we discuss strictly stable case. 
Section \ref{sec: A} is devoted to the study of the asymmetric case.

\section{General discussions} \label{sec: G}

Let $ \{ X(t),\bP_x \} $ be a L\'evy process on $ \bR $. 
We have $ \bE_0[\e^{i \lambda X(t)}] = \e^{- t \Psi(\lambda)} $ for $ \lambda \in \bR $, 
where the L\'evy--Khintchine exponent $ \Psi(\lambda) $ may be represented as 
\begin{align}
\Psi(\lambda) = ib \lambda + v \lambda^2 
- \int_{\bR} \rbra{ \e^{i \lambda x} - 1 - i \lambda x 1_{(-1,1)}(x) } \nu(\d x) 
\label{}
\end{align}
for some constants $ b \in \bR $ and $ v \ge 0 $ 
and some measure $ \nu $ on $ \bR $ such that 
$ \nu(\{ 0 \})=0 $ and $ \int_{\bR} (|x|^2 \wedge 1) \nu(\d x) < \infty $. 
The real and imaginary parts are denoted by $ \theta(\lambda) $ and $ \omega(\lambda) $, 
respectively: 
\begin{align}
\theta(\lambda) =& \Re \Psi(\lambda) 
= v \lambda^2 + \int_{\bR} \rbra{ 1 - \cos \lambda x } \nu(\d x) , 
\label{} \\
\omega(\lambda) =& \Im \Psi(\lambda) 
= b \lambda + \int_{\bR} \rbra{ \lambda x 1_{(-1,1)}(x) - \sin \lambda x } \nu(\d x) . 
\label{}
\end{align}
Note that $ \theta(\lambda) \ge 0 $, 
$ \theta(-\lambda) = \theta(\lambda) $ 
and $ \omega(-\lambda) = - \omega(\lambda) $ 
for all $ \lambda \in \bR $. 

For $ x \in \bR $, we denote 
\begin{align}
T_x = \inf \{ t>0 : X(t)=x \}. 
\label{}
\end{align}
We say that {\em the state 0 is regular for itself} if 
$ \bP_0(T_0=0)=1 $. 
The following theorem is due to Kesten \cite{MR0272059} and Bretagnolle \cite{MR0368175}. 

\begin{Thm}[{\cite{MR0272059} and \cite{MR0368175}}]
The process $ \{ X(t),\bP_x \} $ is not a compound Poisson process 
and is such that the state 0 is regular for itself 
if and only if the following conditions hold: 
\begin{description}
\item[(L1)] 
$ \displaystyle \int_{\bR} \Re \frac{1}{q + \Psi(\lambda)} \d \lambda < \infty $ 
for all $ q>0 $; 
\item[(L2)] 
either $ v>0 $ or $ \displaystyle \int_{(-1,1)} |x| \nu(\d x) = \infty $. 
\end{description}
\end{Thm}

We introduce the following condition which is stronger than {\bf (L1)}: 
\begin{description}
\item[(L1$ ' $)] 
$ \displaystyle \int_0^{\infty } \frac{1}{q+\theta(\lambda)} \d \lambda < \infty $ 
for all $ q>0 $. 
\end{description}
By {\bf (L1$ ' $)}, we have 
\begin{align}
|\e^{-t \Psi(\lambda)}| = \e^{-t \theta(\lambda)} \in L^1(\d \lambda) , 
\label{}
\end{align}
and hence we may invert the Fourier transform so that the function 
\begin{align}
p_t(x) = \frac{1}{2 \pi} \int_{\bR} \Re \e^{- i \lambda x - t \Psi(\lambda)} \d \lambda 
, \quad t>0 , \ x \in \bR 
\label{}
\end{align}
is a continuous transition density with respect to the Lebesgue measure, i.e., 
\begin{align}
\bE_x[f(X(t))] = \int_{\bR} f(y) p_t(y-x) \d y . 
\label{}
\end{align}
Again by {\bf (L1$ ' $)}, 
the Laplace transform of $ p_t(x) $ can be expressed as 
\begin{align}
r_q(x) 
:=& \int_0^{\infty } \e^{-qt} p_t(x) \d t 
\label{} \\
=& \frac{1}{2 \pi} \int_{\bR} \Re \frac{\e^{-i \lambda x}}{q+\Psi(\lambda)} \d \lambda 
\label{} \\
=& \frac{1}{2 \pi} \int_{\bR} \frac{(q+\theta(\lambda)) \cos \lambda x + \omega(\lambda) \sin \lambda x }{(q+\theta(\lambda))^2 + \omega(\lambda)^2} \d \lambda 
\label{}
\end{align}
and thus $ r_q(x) $ is a continuous resolvent density, i.e., 
\begin{align}
R_qf(x) 
:= \bE_x \sbra{ \int_0^{\infty } \e^{-qt} f(X(t)) \d t } 
= \int_{\bR} f(y) r_q(y-x) \d y . 
\label{}
\end{align}

We now recall the excursion theory. 
Let $ \bD $ denote the set of c\`adl\`ag paths on $ [0,\infty ) $ 
taking values in $ \bR $. 
We set 
\begin{align}
\bD^0 = \cbra{ w \in \bD : w(t) = w(t \wedge T_0) \ \text{for all $ t \ge 0 $}} , 
\label{}
\end{align}
which is sometimes called the {\em excursion space}. 
We denote the coordinate process on $ \bD $ by $ X(t) $: 
$ X(t)(w) = w(t) $ for $ w \in \bD $. 
The set $ \bD $ is equipped with the $ \sigma $-field generated by 
the totality of cylinder sets of $ \bD $. 
We denote 
\begin{align}
\cF_t = \sigma(X(s): s \le t) . 
\label{}
\end{align}

Since the state 0 is regular for itself, there exists a local time at 0, 
i.e., a positive continuous additive functional 
which increases only when $ X(t) $ stays at $ 0 $. 
We denote by $ L(t) $ the version of the local time normalized so that 
\begin{align}
\bE_x \sbra{ \int_0^{\infty } \e^{-qt} \d L(t) } = r_q(-x) 
\quad \text{for all $ q>0 $ and $ x \in \bR $}. 
\label{}
\end{align}
We write $ \eta(l) $ for the right-continuous inverse of the local time: 
\begin{align}
\eta(l) = \inf \{ t \ge 0 : L(t) > l \} 
\quad \text{for $ l \ge 0 $}. 
\label{}
\end{align}
Then the process $ \{ \eta(l),\bP_0 \} $ is a subordinator 
which explodes at time $ \lambda := L(\infty ) $ when $ X(t) $ is transient. 

We now recall the excursion theory for our L\'evy process. 
Let $ \vn $ be a $ \sigma $-finite measure on $ \bD $ such that 
\begin{align}
\vn[T_0 \wedge 1] < \infty 
\quad \text{and} \quad 
\vn(\bD) = \infty . 
\label{eq: cond vn}
\end{align}
Let $ p(l)=(p(l;t))_{t \ge 0} $ be a Poisson point process with characteristic measure $ \vn $ 
defined on a probability space $ (\Omega,\cF,\bP) $; 
see \cite{MR0402949} for the precise definition. 
We set 
\begin{align}
\eta_p(l) = \sum_{s \le l} T_0(p(s)) 
\label{}
\end{align}
and set 
\begin{align}
\lambda_p = \inf \{ l>0 : \eta_p(l)=\infty \} . 
\label{}
\end{align}
By the condition \eqref{eq: cond vn}, we see that 
$ \eta_p(l) $ is strictly-increasing in $ l \in [0,\lambda_p) $, 
and hence that there exists a continuous inverse process $ L_p(t) $. 
We define 
\begin{align}
X_p(t) = 
\begin{cases}
p(l;t-\eta_p(l-)) 
& \text{if $ \eta_p(l-) \le t < \eta_p(l) $} \\
& \ \ \text{for some $ l \in [0,\lambda_p] $}, \\
0 & \text{otherwise}. 
\end{cases}
\label{}
\end{align}
The following theorem is due to It\^o \cite{MR0402949} and Meyer \cite{MR0388551}. 

\begin{Thm}[{\cite{MR0402949} and \cite{MR0388551}}]
There exists a unique $ \sigma $-finite measure on $ \bD^0 $ such that 
the condition \eqref{eq: cond vn} is satisfied 
and that, for a Poisson point process $ p(l) $ with characteristic measure $ \vn $ 
defined on a probability space $ (\Omega,\cF,\bP) $, 
it follows that 
\begin{align}
\{ \eta_p, L_p, X_p, \lambda_p, \bP \} 
\law 
\{ \eta, L, X, \lambda, \bP_0 \} . 
\label{}
\end{align}
\end{Thm}

Using compensation formula, we can obtain the following. 

\begin{Thm} \label{thm: ent law}
For any non-negative measurable function $ f $, one has 
\begin{align}
\vn \sbra{ \int_0^{T_0} \e^{-qt} f(X(t)) \d t } 
= \frac{1}{r_q(0)} \int_{\bR} f(y) r_q(y) \d y . 
\label{}
\end{align}
Consequently, taking $ f \equiv 1 $, one has 
\begin{align}
\int_0^{\infty } \e^{-qt} \vn(T_0>t) \d t 
= \frac{1}{q r_q(0)} . 
\label{}
\end{align}
\end{Thm}

A non-negative Borel function $ h $ on $ \bR $ is called 
{\em harmonic for the killed process} if the following conditions hold: 
\begin{enumerate}
\item 
$ \bE^0_x[h(X(t))] = h(x) $ for all $ x \in \bR $; 
\item 
$ \vn[h(X(t))] $ is finite constant in $ t>0 $. 
\end{enumerate}
It is easy to see that $ h $ is harmonic for the killed process if and only if 
\begin{align}
qR_qh(x) = h(x) + r_q(-x) 
\quad \text{for all $ q>0 $ and $ x \in \bR $}. 
\label{eq: harm func}
\end{align}
If $ h $ is harmonic for the killed process 
and if $ h(x)>0 $ for all $ x \neq 0 $, then 
we may define {\em $ h $-transform} $ \{ X(t),\bP^h_x \} $ by 
\begin{align}
\bE^h_x \sbra{ F_t } = 
\begin{cases}
\displaystyle 
\frac{1}{h(x)} \bE^0_x \sbra{ F_t h(X(t)) } 
& \text{if $ x \neq 0 $}, \\
\vn \sbra{ F_t h(X(t)) } 
& \text{if $ x=0 $} 
\end{cases}
\label{}
\end{align}
for all $ t>0 $ and all non-negative $ \cF_t $-measurable functional $ F_t $. 
Note that the process $ \{ X(t),\bP^h_x \} $ may be regarded 
as the process conditioned to avoid zero; see \cite{MR2648275} 
for the precise explanation. 
The following theorem is a generalization of the formula \eqref{eq: BM disinteg}. 

\begin{Thm} \label{thm: disinteg}
Suppose that the following conditions are satisfied: 
\begin{description}
\item[(H)] 
$ h $ is harmonic for the killed process and 
$ h(x)>0 $ for all $ x \neq 0 $; 
\item[(T)] 
there exists a density $ \rho(t) $ such that $ \vn(T_0 \in \d t) = \rho(t) \d t $ 
on $ (0,\infty ) $. 
\end{description}
Then the excursion measure $ \vn $ enjoys the following disintegration formula: 
\begin{align}
\vn(\cdot) = \int_0^{\infty } \bP^h_0(\cdot | X(t-)=0) \rho(t) \d t 
+ \kappa \bP^h_0(\cdot) , 
\label{}
\end{align}
where 
\begin{align}
\kappa = \vn(T_0=\infty ) = \lim_{q \to 0+} \frac{1}{r_q(0)} . 
\label{}
\end{align}
\end{Thm}

\section{Renormalized zero resolvent} \label{sec: R}

For $ q>0 $, we set 
\begin{align}
h_q(x) = r_q(0)-r_q(-x) 
, \quad x \in \bR . 
\label{}
\end{align}
We call $ h_0 := \lim_{q \to 0+} h_q $ the {\em renormalized zero resolvent}, 
if the limit exists. 
We shall prove the following proposition. 

\begin{Prop} \label{prop: h}
Suppose that $ h_0 $ exists. Then 
\begin{align}
R_qh_0(x) < \infty 
\quad \text{for all $ q>0 $ and all $ x \in \bR $}. 
\label{}
\end{align}
Suppose, moreover, that there exists a non-negative measurable function $ f(x) $ such that 
\begin{align}
R_qf(x)<\infty 
\quad \text{for all $ q>0 $ and all $ x \in \bR $} 
\label{}
\end{align}
and 
\begin{align}
h_q(x) \le f(x) 
\quad \text{for all $ q>0 $ and all $ x \in \bR $}. 
\label{eq: h ass}
\end{align}
Then $ h_0 $ is harmonic for the killed process. 
\end{Prop}

To prove Proposition \ref{prop: h}, we need the following lemma. 

\begin{Lem} \label{lem: theta}
Suppose that the conditions {\bf (L1$ ' $)} and {\bf (L2)} are satisfied. 
Then the following assertions hold: 
\begin{enumerate}
\item 
$ \theta(\lambda) \to \infty $ 
as $ \lambda \to \infty \Big. $. 
\item 
$ \displaystyle \int_0^1 \frac{\lambda^2}{\theta(\lambda)} \d \lambda 
+ \int_1^{\infty } \frac{1}{\theta(\lambda)} \d \lambda < \infty $. 
\item 
$ \displaystyle \lim_{q \to 0+} q r_q(0) = 0 \Big. $. 
\end{enumerate}
\end{Lem}

\Proof{
(i) 
Since $ \e^{- \Psi(\lambda)} = \int_{-\infty }^{\infty } \e^{i \lambda x} p_1(x) \d x $, 
the Riemann--Lebesgue theorem yields 
$ \e^{- \theta(\lambda)} = |\e^{- \Psi(\lambda)}| \to 0 $ 
as $ \lambda \to \infty $. 
This proves Claim (i). 

(ii) By the dominated convergence theorem, we have 
\begin{align}
\frac{\theta(\lambda)}{\lambda^2} 
\ge& v + \int_{(-1,1)} \frac{1-\cos \lambda x}{(\lambda x)^2} x^2 \nu(\d x) 
\label{} \\
\tend{}{\lambda \to 0}& 
v + \frac{1}{2} \int_{(-1,1)} x^2 \nu(\d x) > 0 , 
\label{}
\end{align}
which shows $ \int_0^1 \frac{\lambda^2}{\theta(\lambda)} \d \lambda < \infty $. 
By (i), we have $ \theta(\lambda)/(1+\theta(\lambda)) \to 1 $ 
as $ \lambda \to \infty $. 
Hence, by {\bf (L1$ ' $)}, we obtain 
$ \int_1^{\infty } \frac{1}{\theta(\lambda)} \d \lambda < \infty $. 

(iii) 
Since we have 
\begin{align}
\absol{ \Re \frac{q}{q+\Psi(\lambda)} } 
\le \frac{q}{q+\theta(\lambda)} 
\le \min \cbra{ 1 , \frac{1}{\theta(\lambda)} } , 
\label{}
\end{align}
we may apply the dominated convergence theorem and we obtain 
\begin{align}
\lim_{q \to 0+} 2 \pi q r_q(0) 
=& \lim_{q \to 0+} \int_{\bR} \Re \frac{q}{q+\Psi(\lambda)} \d \lambda 
\label{} \\
=& \int_{\bR} \lim_{q \to 0+} \Re \frac{q}{q+\Psi(\lambda)} \d \lambda 
\label{} \\
=& 0 . 
\label{}
\end{align}
The proof is now complete. 
}

We now prove Proposition \ref{prop: h}. 

\Proof[Proof of Proposition \ref{prop: h}]{
Let $ q>0 $ and $ p>0 $. 
Using the resolvent equation 
\begin{align}
r_q(z-x) - r_p(z-x) + (q-p) \int_{\bR} r_q(y-x) r_p(z-y) \d y = 0 
\label{}
\end{align}
for $ x,z \in \bR $, we have 
\begin{align}
R_qh_p(x) 
= \frac{h_p(x)}{q-p} + \frac{r_q(-x)}{q-p} - \frac{pr_p(0)}{q(q-p)} . 
\label{eq: hp id}
\end{align}
By Fatou's lemma and by (iii) of Lemma \ref{lem: theta}, we have 
\begin{align}
q R_qh_0(x) 
\le \liminf_{p \to 0+} q R_qh_p(x) 
= h_0(x) + r_q(-x) < \infty . 
\label{eq: h integ}
\end{align}

By the assumption \eqref{eq: h ass} and by the integrability \eqref{eq: h integ}, 
we may apply the dominated convergence theorem 
to \eqref{eq: hp id} being taken limit as $ q \to 0+ $, 
and thus we obtain the identity \eqref{eq: harm func} for $ h=h_0 $. 
This proves the claim. 
}

\section{Symmetric case} \label{sec: S}

Let us recall several results of the symmetric case obtained in \cite{MR2603019}. 
Assume that $ \{ X(t),\bP_x \} $ is symmetric, 
i.e., $ \{ X(t),\bP_x \} $ is identical in law to $ \{ -X(t),\bP_{-x} \} $. 
In this case, we have $ \nu(-B)=\nu(B) $ for all $ B \in \cB(\bR) $, 
where $ -B = \{ -x \in \bR : x \in B \} $. 
Hence we have $ \omega(\lambda) \equiv 0 $ 
and 
\begin{align}
r_q(x) = \frac{1}{\pi} \int_0^{\infty } \frac{\cos \lambda x}{q+\theta(\lambda)} \d \lambda 
\quad \text{for all $ q>0 $ and $ x \in \bR $}. 
\label{}
\end{align}

\begin{Thm}[{\cite[Theorem 3.2]{MR2603019}}] \label{thm: sym Levy rho}
Suppose that $ \{ X(t),\bP_x \} $ is symmetric 
and that the conditions {\bf (L1$ ' $)} and {\bf (L2)} are satisfied. 
Then there exists a $ \sigma $-finite measure $ \sigma^* $ on $ [0,\infty ) $ such that 
\begin{align}
\frac{1}{q r_q(0)} = \int_{[0,\infty )} \frac{1}{q+\xi} \sigma^*(\d \xi) 
\quad \text{for all $ q>0 $}. 
\label{eq: sym Levy qrq}
\end{align}
Consequently, it holds that $ \kappa = \sigma^*(\{ 0 \}) $ 
and the condition {\bf (T)} is satisfied with 
\begin{align}
\rho(t) = \int_{(0,\infty )} \e^{-t \xi} \xi \sigma^*(\d \xi) . 
\label{}
\end{align}
\end{Thm}

\begin{Thm}[{\cite[Theorems 1.1 and 1.2]{MR2603019}}] \label{thm: sym Levy h0}
Suppose that $ \{ X(t),\bP_x \} $ is symmetric 
and that the conditions {\bf (L1$ ' $)} and {\bf (L2)} are satisfied. 
Then $ h_0 $ exists and is given as 
\begin{align}
h_0(x) = \frac{1}{\pi} \int_0^{\infty } \frac{1-\cos \lambda x}{\theta(\lambda)} \d \lambda 
\quad \text{for all $ x \in \bR $}. 
\label{eq: sym Levy h0}
\end{align}
Moreover, $ h=h_0(x) $ is harmonic for the killed process 
and the condition {\bf (H)} is satisfied. 
\end{Thm}

Here we note that Theorem \ref{thm: sym Levy h0} is immediate 
from Proposition \ref{prop: h} and the fact that $ h_q(x) $ increases as $ q $ decreases to 0, 
which can be seen by the expression 
\begin{align}
h_q(x) 
= \frac{1}{\pi} \int_0^{\infty } \frac{1-\cos \lambda x}{q + \theta(\lambda)} \d \lambda . 
\label{}
\end{align}
The dominating function $ f $ in Proposition \ref{prop: h} 
can be chosen to be $ h_0 $.

\section{Strictly stable case} \label{sec: SS} 

Let us consider the case where 
$ \{ X(t),\bP_x \} $ is a strictly stable process and is not a Brownian motion. 
The L\'evy measure $ \nu $ is given as 
\begin{align}
\nu(\d x) = 
\begin{cases}
c_+ |x|^{-\alpha -1} \d x 
\quad \text{on $ (0,\infty ) $}, \\
c_- |x|^{-\alpha -1} \d x 
\quad \text{on $ (-\infty ,0) $}, 
\end{cases}
\label{}
\end{align}
where $ 1<\alpha <2 $ 
and $ c_+ $ and $ c_- $ are non-negative constants 
such that $ c_++c_->0 $, 
and the L\'evy--Khintchine exponent is given as 
\begin{align}
\Psi(\lambda) = c_{\theta} |\lambda|^{\alpha } 
\rbra{ 1 - i \beta \sgn(\lambda) \tan \frac{\pi \alpha }{2} } , 
\label{}
\end{align}
where $ c_{\theta}>0 $ and $ \beta \in [-1,1] $ are constants given as 
\begin{align}
c_{\theta} = (c_++c_-) \frac{\pi}{\alpha S_{\alpha }} 
, \quad 
\beta = \frac{c_+-c_-}{c_++c_-} , 
\label{}
\end{align}
with $ S_{\alpha } = 2 \Gamma(\alpha ) \sin \frac{\pi \alpha }{2} $. 
We then have 
\begin{align}
\theta(\lambda) = c_{\theta} |\lambda|^{\alpha } 
\quad \text{and} \quad 
\omega(\lambda) = c_{\omega} |\lambda|^{\alpha } \sgn(\lambda) , 
\label{eq: stable theta and omega}
\end{align}
where 
$ c_{\omega} = c_{\theta} \beta \cdot (- \tan \frac{\pi \alpha }{2}) $. 

By \eqref{eq: stable theta and omega}, we have 
\begin{align}
p_t(0) = c_p t^{-\frac{1}{\alpha }} 
, \quad 
r_q(0) = c_r q^{\frac{1}{\alpha }-1} , 
\label{}
\end{align}
where $ c_p $ and $ c_r $ are constants given as 
\begin{align}
c_p =& \frac{1}{\pi} \int_0^{\infty } \rbra{ \cos c_{\omega} \lambda^{\alpha } } \e^{- c_{\theta} \lambda^{\alpha }} \d \lambda 
, \label{} \\
c_r =& \frac{1}{\pi} \int_0^{\infty } \frac{1+c_{\theta} \lambda^{\alpha }}{(1+c_{\theta} \lambda^{\alpha })^2 + (c_{\omega} \lambda^{\alpha })^2} \d \lambda . 
\label{}
\end{align}
The condition {\bf (T)} is satisfied with $ \rho(t) $ being given as 
\begin{align}
\rho(t) = \frac{1-\frac{1}{\alpha }}{c_r \Gamma(\frac{1}{\alpha })} t^{\frac{1}{\alpha }-2} . 
\label{}
\end{align}
By the direct computation, we can obtain the following formula: 
\begin{align}
h_0(x) = \frac{1 - \beta \sgn(x)}{c_{\theta} \rbra{1 + \beta^2 \tan^2 \frac{\pi \alpha }{2}} C_{\alpha } } |x|^{\alpha -1} , 
\label{eq: Port}
\end{align}
where $ C_{\alpha } = 2 \Gamma(\alpha ) \cdot (-\cos \frac{\pi \alpha }{2}) $. 

Note that this function $ h_0 $ has been appeared in Port \cite[Theorem 2]{MR0217877}: 
in the case $ -1<\beta<1 $, he obtained the asymptotic result 
\begin{align}
\bP_x(T_0>t) \sim C(\alpha ,\beta) t^{- \rho} h_0(x) 
\quad \text{as $ t \to \infty $}, 
\label{}
\end{align}
where $ \rho = 1-\frac{1}{\alpha } \in (\frac{1}{2},1) $ 
and $ C(\alpha ,\beta) $ is a certain positive constant 
and where $ f \sim g $ means $ f/g \to 1 $. 
In the case $ \beta=1 $, i.e., 
the spectrally positive case, we have 
\begin{align}
h_0(x) = 
\begin{cases}
0 & \text{if $ x \ge 0 $}, \\
c(\alpha) |x|^{\alpha -1} & \text{if $ x<0 $} 
\end{cases}
\label{}
\end{align}
for some constant $ c(\alpha )>0 $. 
In this case, 
we consult Port \cite[Theorem 2]{MR0217877} for the case $ x<0 $ 
and Chaumont \cite[third line on pp.386]{MR1465814} for the case $ x>0 $ 
and we obtain 
\begin{align}
\bP_x(T_0>t) \sim 
\begin{cases}
C_+(\alpha ,1) t^{-\frac{1}{\alpha }} x  & \text{if $ x>0 $}, \\
C_-(\alpha ,1) t^{-\rho} |x|^{\alpha -1} & \text{if $ x<0 $} 
\end{cases}
\label{}
\end{align}
as $ t \to \infty $, 
where $ C_+(\alpha ,1) $ and $ C_-(\alpha ,1) $ are some positive constants. 

For harmonicity of $ h_0 $, we obtain the following. 

\begin{Thm} \label{thm: SS1} 
Suppose that $ \{ X(t),\bP_x \} $ is a strictly stable process of index $ 1<\alpha <2 $. 
If $ -1<\beta<1 $, then the condition {\bf (H)} is satisfied. 
If $ \beta = 1 $ or $ -1 $, 
then $ h_0 $ is harmonic for the killed process. 
\end{Thm}

Theorem \ref{thm: SS1} will be proved in Section \ref{sec: H}. 

\begin{Rem}
In the case $ \beta=1 $, Silverstein proved in \cite[Theorem 2]{MR573292} that 
the function $ h_0 $ satisfies 
\begin{align}
\bE_x \sbra{ h_0(X(t)) ; T_{[0,\infty )} > t } 
= h_0(x) 
\quad \text{for all $ x<0 $}, 
\label{}
\end{align}
where $ T_{[0,\infty )} = \inf \{ t>0 : X(t) \in [0,\infty ) \} $. 
We can recover this result from Theorem \ref{thm: SS1}. 
In fact, 
once the process starting from $ x<0 $ changes sign, 
then it stays positive until hitting zero. 
We thus obtain 
\begin{align}
\bE_x \sbra{ h_0(X(t)) ; T_{[0,\infty )} > t } 
= \bE^0_x \sbra{ h_0(X(t)) } 
= h_0(x) 
\label{}
\end{align}
for all $ x<0 $. We also note that 
the transformed process $ \{ X(t),\bP^{h_0}_x \} $, which is defined only for $ x<0 $, 
is nothing else but the process conditioned to stay negative, 
which was introduced by Chaumont \cite{MR1419491}. 
\end{Rem}

\section{Asymmetric case} \label{sec: A}

\subsection{Existence of renormalized zero resolvent}

\begin{Prop} \label{prop: asymp Levy}
Suppose that the conditions {\bf (L1$ ' $)} and {\bf (L2)} are satisfied. 
Suppose, in addition, that the following condition is satisfied: 
\begin{description}
\item[(L3)] 
$ \theta $ and $ \omega $ have measurable derivatives on $ (0,\infty ) $ which satisfy 
\begin{align}
\int_0^{\infty } 
\frac{(|\theta'(\lambda)| + |\omega'(\lambda)|) (\lambda^2 \wedge 1)}
{\theta(\lambda)^2 + \omega(\lambda)^2} 
\d \lambda < \infty . 
\label{}
\end{align}
\end{description}
Then $ h_0 $ exists and is given as 
\begin{align}
h_0(x) 
= \frac{1}{\pi} \int_0^{\infty } \frac{\theta(\lambda) (1-\cos \lambda x) + \omega(\lambda) \sin \lambda x }{\theta(\lambda)^2 + \omega(\lambda)^2} \d \lambda . 
\label{}
\end{align}
\end{Prop}

\Proof{
Note that 
\begin{align}
h_q(x) 
= r_q(0) - r_q(-x) 
= h^{(1)}_q(x) + h^{(2)}_q(x) , 
\label{}
\end{align}
where 
\begin{align}
h^{(1)}_q(x) 
=& \frac{1}{\pi} \int_0^{\infty } \frac{(q+\theta(\lambda)) (1-\cos \lambda x) }{(q+\theta(\lambda))^2 + \omega(\lambda)^2} \d \lambda , 
\label{} \\
h^{(2)}_q(x) 
=& \frac{1}{\pi} \int_0^{\infty } \frac{\omega(\lambda) \sin \lambda x }{(q+\theta(\lambda))^2 + \omega(\lambda)^2} \d \lambda . 
\label{}
\end{align}

By (ii) of Lemma \ref{lem: theta}, we have 
\begin{align}
\int_0^{\infty } \frac{1-\cos \lambda x}{\theta(\lambda)} \d \lambda < \infty 
\quad \text{for all $ x \in \bR $}. 
\label{}
\end{align}
Hence we may apply the dominated convergence theorem and see that 
$ h^{(1)}_q(x) \to h^{(1)}_0(x) $ as $ q \to 0+ $, where 
\begin{align}
h^{(1)}_0(x) =& \frac{1}{\pi} 
\int_0^{\infty } \frac{\theta(\lambda) (1-\cos \lambda x) }{\theta(\lambda)^2 + \omega(\lambda)^2} \d \lambda . 
\label{}
\end{align}

Suppose that $ x \neq 0 $. 
Integrating by parts, we have 
\begin{align}
h^{(2)}_q(x) 
= \frac{1}{\pi x} \int_0^{\infty } \frac{f_q(\lambda) (1-\cos \lambda x) }{\{ (q+\theta(\lambda))^2 + \omega(\lambda)^2 \}^2} \d \lambda , 
\label{}
\end{align}
where 
\begin{align}
\begin{split}
f_q(\lambda) 
=& \omega(\lambda) \frac{\partial }{\partial \lambda} \{ (q+\theta(\lambda))^2 + \omega(\lambda)^2 \} 
\\
&- \omega'(\lambda) \{ (q+\theta(\lambda))^2 + \omega(\lambda)^2 \} 
\end{split}
\label{} \\
\begin{split}
=& \omega(\lambda) \cdot 2 (q+\theta(\lambda)) \theta'(\lambda) 
\\
&- \omega'(\lambda) \{ (q+\theta(\lambda))^2 + \omega(\lambda)^2 \} . 
\end{split}
\label{}
\end{align}
We now have 
\begin{align}
|f_q(\lambda)| 
\le (|\theta'(\lambda)| + |\omega'(\lambda)|) 
\{ (q+\theta(\lambda))^2 + \omega(\lambda)^2 \} 
\label{}
\end{align}
and hence we have 
\begin{align}
\frac{|f_q(\lambda)| (1-\cos \lambda x) }{\{ (q+\theta(\lambda))^2 + \omega(\lambda)^2 \}^2} 
\le 
\frac{(|\theta'(\lambda)| + |\omega'(\lambda)|) (1-\cos \lambda x) }
{ \theta(\lambda)^2 + \omega(\lambda)^2 } . 
\label{eq: fq}
\end{align}
Therefore, by the dominated convergence theorem, we obtain 
$ h^{(2)}_q(x) \to h^{(2)}_0(x) $ as $ q \to 0+ $, where 
\begin{align}
h^{(2)}_0(x) 
= \frac{1}{\pi x} \int_0^{\infty } \frac{f_0(\lambda) (1-\cos \lambda x) }{\{ \theta(\lambda)^2 + \omega(\lambda)^2 \}^2} \d \lambda . 
\label{}
\end{align}
Again by integration by parts, we obtain the desired result. 
}

\subsection{Harmonicity of renormalized zero resolvent} \label{sec: H}

\begin{Thm} \label{thm: AL}
Suppose that $ \theta $ and $ \omega $ have measurable derivatives on $ (0,\infty ) $ 
which satisfy the following conditions: 
\begin{enumerate}
\item 
there exist constants $ \un{c}_{\theta}>0 $ and $ \bar{c}_{\theta}>0 $ such that 
\begin{align}
\alpha \un{c}_{\theta} \lambda^{\alpha -1} 
\le \theta'(\lambda) \le 
\alpha \bar{c}_{\theta} \lambda^{\alpha -1} 
\quad \text{for all $ \lambda>0 $}; 
\label{}
\end{align}
\item 
there exist constants $ \un{c}_{\omega}>0 $ and $ \bar{c}_{\omega}>0 $ such that 
\begin{align}
\alpha \un{c}_{\omega} \lambda^{\alpha -1} 
\le \omega'(\lambda) \le 
\alpha \bar{c}_{\omega} \lambda^{\alpha -1} 
\quad \text{for all $ \lambda>0 $}. 
\label{}
\end{align}
\end{enumerate}
Then $ h_0 $ exists and there exists a constant $ C $ such that 
\begin{align}
h_q(x) \le C |x|^{\alpha -1} 
\quad \text{for all $ q>0 $ and all $ x \in \bR $}. 
\label{}
\end{align}
Suppose, moreover, that the following condition is satisfied: 
\begin{enumerate}
\setcounter{enumi}{2}
\item 
it holds that 
\begin{align}
\frac{\un{c}_{\theta} (\un{c}_{\theta}^2 + \un{c}_{\omega}^2) }{ \bar{c}_{\theta}^2 + \bar{c}_{\omega}^2 } \cdot (- \tan \frac{\pi \alpha }{2}) 
> \max \{ \bar{c}_{\omega} , \bar{c}_{\theta} -\un{c}_{\omega} \} . 
\label{}
\end{align}
\end{enumerate}
Then one has 
\begin{align}
\un{c} |x|^{\alpha -1} \le h_0(x) \le \bar{c} |x|^{\alpha -1} 
, \quad x \in \bR 
\label{eq: ineq asym Levy}
\end{align}
for some positive constants $ \un{c} $ and $ \bar{c} $. 
Consequently, $ h_0 $ is harmonic for the killed process 
and the condition {\bf (H)} is satisfied. 
\end{Thm}

\Proof{
We recall the following formula (see, e.g., \cite[Proposition 7.1]{MR2599211}): 
\begin{align}
C_{\alpha } 
:=& \int_0^{\infty } \frac{1-\cos x}{x^{\alpha }} \d x 
\label{} \\
=& \frac{\pi}{2 \Gamma(\alpha ) (-\cos \frac{\pi \alpha }{2})} 
\quad \text{for $ \alpha \in (1,3) $}. 
\label{}
\end{align}
Note that 
\begin{align}
\int_0^{\infty } \frac{x-\sin x}{x^{\alpha +1}} \d x 
= \frac{C_{\alpha }}{\alpha } 
\quad \text{for $ \alpha \in (1,3) $} 
\label{}
\end{align}
and 
\begin{align}
\frac{C_{\alpha }}{\alpha C_{\alpha +1}} 
= - \tan \frac{\pi \alpha }{2} 
\quad \text{for $ \alpha \in (1,2) $}. 
\label{}
\end{align}

Suppose that the conditions (i) and (ii) are satisfied. 
Since we have $ \theta(0) = \omega(0) = 0 $, we have 
\begin{align}
\un{c}_{\theta} \lambda^{\alpha } 
\le \theta(\lambda) \le 
\bar{c}_{\theta} \lambda^{\alpha } 
\quad \text{for all $ \lambda \ge 0 $} 
\label{}
\end{align}
and 
\begin{align}
\un{c}_{\omega} \lambda^{\alpha } 
\le \omega(\lambda) \le 
\bar{c}_{\omega} \lambda^{\alpha } 
\quad \text{for all $ \lambda \ge 0 $}. 
\label{}
\end{align}
We now notice that 
the conditions {\bf (L1$ ' $)}, {\bf (L2)} and {\bf (L3)} are all satisfied, 
which yields by Proposition \ref{prop: asymp Levy} that $ h_0 $ exists. 
From the proof of Proposition \ref{prop: asymp Levy}, it holds, 
for all $ q>0 $ and $ x \neq 0 $, that 
\begin{align}
h_q(x) 
\le& \frac{1}{\pi} \int_0^{\infty } \frac{1-\cos \lambda x}{\theta(\lambda)} \d \lambda 
+ \frac{1}{\pi |x|} \int_0^{\infty } 
\frac{4c_3(1-\cos \lambda x)}{\lambda \theta(\lambda)} \d \lambda 
\label{} \\
\le& \frac{1}{\pi} \cdot \frac{C_{\alpha }}{\un{c}_{\theta}} |x|^{\alpha -1} 
+ \frac{4c_3}{\pi |x|} \cdot \frac{C_{\alpha +1}}{\un{c}_{\theta}} |x|^{\alpha } 
\label{} \\
=& c_5 |x|^{\alpha -1} 
\label{}
\end{align}
for some constant $ c_5>0 $. 

Suppose that the condition (iii) is also satisfied. 
To prove the claim, 
it suffices to prove that $ h_0(x) \ge \un{c} |x|^{\alpha -1} $ 
for some constant $ \un{c}>0 $. 
Note that 
\begin{align}
h^{(1)}_0(x) 
\ge& \frac{\un{c}_{\theta} }{\pi (\bar{c}_{\theta}^2 + \bar{c}_{\omega}^2)} 
\int_0^{\infty } \frac{1-\cos \lambda x}{ \lambda^{\alpha } } \d \lambda 
\label{} \\
=& \frac{\un{c}_{\theta} }{\pi (\bar{c}_{\theta}^2 + \bar{c}_{\omega}^2)} 
C_{\alpha } |x|^{\alpha -1} . 
\label{}
\end{align}
First, we let $ x>0 $. Since we have 
\begin{align}
- f_0(\lambda) 
=& 2 (- \omega(\lambda)) \theta(\lambda) \theta'(\lambda) 
+ \{ \theta(\lambda)^2 + \omega(\lambda)^2 \} \omega'(\lambda) 
\label{} \\
\le& \{ \theta(\lambda)^2 + \omega(\lambda)^2 \} 
\cdot \alpha \bar{c}_{\omega} \lambda^{\alpha -1} , 
\label{}
\end{align}
we have 
\begin{align}
- h^{(2)}_0(x) 
\le& \frac{1}{\pi x} \int_0^{\infty } \frac{\alpha \bar{c}_{\omega} \lambda^{\alpha -1} (1-\cos \lambda x)}{\theta(\lambda)^2 + \omega(\lambda)^2} \d \lambda 
\label{} \\
\le& \frac{1}{\pi x} \cdot \frac{\alpha \bar{c}_{\omega} }{\un{c}_{\theta}^2 + \un{c}_{\omega}^2} \int_0^{\infty } \frac{1-\cos \lambda x}{\lambda^{\alpha +1}} \d \lambda 
\label{} \\
=& \frac{\bar{c}_{\omega} }{ \pi (\un{c}_{\theta}^2 + \un{c}_{\omega}^2) } \alpha C_{\alpha +1} |x|^{\alpha -1} . 
\label{}
\end{align}
Hence we have 
\begin{align}
h_0(x) \ge 
\rbra{ \frac{\un{c}_{\theta} }{ \bar{c}_{\theta}^2 + \bar{c}_{\omega}^2 } 
(- \tan \frac{\pi \alpha }{2}) 
- \frac{\bar{c}_{\omega} }{ \un{c}_{\theta}^2 + \un{c}_{\omega}^2 } } 
\frac{\alpha C_{\alpha +1}}{\pi} |x|^{\alpha -1} . 
\label{}
\end{align}
Second, we let $ x<0 $. Since we have 
\begin{align}
f_0(\lambda) 
=& 2 \omega(\lambda) \theta(\lambda) \theta'(\lambda) 
- \{ \theta(\lambda)^2 + \omega(\lambda)^2 \} \omega'(\lambda) 
\label{} \\
\le& \{ \theta(\lambda)^2 + \omega(\lambda)^2 \} 
\{ \theta'(\lambda) - \omega'(\lambda) \} 
\label{} \\
\le& \{ \theta(\lambda)^2 + \omega(\lambda)^2 \} 
\cdot \alpha (\bar{c}_{\theta}-\un{c}_{\omega}) \lambda^{\alpha -1} , 
\label{}
\end{align}
we have 
\begin{align}
- h^{(2)}_0(x) 
\le& \frac{1}{\pi|x|} \int_0^{\infty } \frac{\alpha (\bar{c}_{\theta}-\un{c}_{\omega}) \lambda^{\alpha -1} (1-\cos \lambda x)}{\theta(\lambda)^2 + \omega(\lambda)^2} \d \lambda 
\label{} \\
\le& \frac{1}{\pi |x|} \cdot \frac{\alpha (\bar{c}_{\theta}-\un{c}_{\omega}) }{\un{c}_{\theta}^2 + \un{c}_{\omega}^2} \int_0^{\infty } \frac{1-\cos \lambda x}{\lambda^{\alpha +1}} \d \lambda 
\label{} \\
=& \frac{\bar{c}_{\theta}-\un{c}_{\omega} }{ \pi (\un{c}_{\theta}^2 + \un{c}_{\omega}^2) } \alpha C_{\alpha +1} |x|^{\alpha -1} . 
\label{}
\end{align}
Hence we have 
\begin{align}
h_0(x) \ge 
\rbra{ \frac{\un{c}_{\theta} }{ \bar{c}_{\theta}^2 + \bar{c}_{\omega}^2 } 
(- \tan \frac{\pi \alpha }{2}) 
- \frac{\bar{c}_{\theta} -\un{c}_{\omega} }{ \un{c}_{\theta}^2 + \un{c}_{\omega}^2 } } 
\frac{\alpha C_{\alpha +1}}{\pi} |x|^{\alpha -1} . 
\label{}
\end{align}

The proof is now complete. 
}

We now proceed to prove Theorem \ref{thm: SS1}. 

\Proof[Proof of Theorem \ref{thm: SS1}]{
The proof in the case $ -1<\beta<1 $ has been already done 
by the former half of Theorem \ref{thm: AL} 
and by the expression \eqref{eq: Port} of $ h_0 $. 

The proof in the case $ \beta=1 $ or $ -1 $ 
is a direct consequence of \cite[Theorem 2]{MR573292}, 
but we give the proof for convenience of the reader. 
We shall prove the claim only for the case $ \beta=1 $. 
Denote 
\begin{align}
(x)_+ = \max \{ x,0 \} 
\quad \text{and} \quad 
(x)_- = \max \{ -x,0 \} . 
\label{}
\end{align}
By Proposition \ref{prop: h}, we see that 
$ R_qf_-(x) < \infty $ for all $ q>0 $ and $ x \in \bR $ 
with $ f_-(x) = (x)_-^{\alpha -1} $. 
It is now sufficient to prove that 
$ f_+(x) := (x)_+^{\alpha -1} $ satisfies 
\begin{align}
R_qf_+(x) < \infty 
\quad \text{for all $ q>0 $ and $ x \in \bR $}; 
\label{eq: xvee0}
\end{align}
in fact, if we assume that \eqref{eq: xvee0} is satisfied, 
then we have $ R_qf(x) < \infty $ for all $ q>0 $ and $ x \in \bR $ 
with $ f(x) = |x|^{\alpha -1} $, 
and hence, by Proposition \ref{prop: h}, 
we conclude that $ h_0 $ is harmonic for the killed process. 

To prove \eqref{eq: xvee0}, 
we recall the following formula from Corollary 2 on pp.94 of \cite{MR854867}: 
\begin{align}
\bP_0(X(1) \ge x) \sim Z(\alpha ,\beta) x^{-\alpha -1} 
\quad \text{as $ x \to \infty $}, 
\label{}
\end{align}
where $ Z(\alpha ,\beta) $ is a certain positive constant. 
Let $ x \in \bR $ be fixed. 
Then we have 
\begin{align}
R_qf_+(x) 
=& \int_0^{\infty } \d t \e^{-qt} \bE_x \sbra{ (X(t))_+^{\alpha -1} } 
\label{} \\
=& \int_0^{\infty } \d t \e^{-qt} \bE_0 \sbra{ (x+t^{1/\alpha } X(1))_+^{\alpha -1} } 
\label{} \\
\le& \int_0^{\infty } \d t \e^{-qt} \bE_0 \sbra{ (|x|+t^{1/\alpha } X(1))_+^{\alpha -1} } . 
\label{}
\end{align}
It is obvious that 
\begin{align}
& \int_0^{\infty } \d t \e^{-qt} 
\bE_0 \sbra{ (|x|+t^{1/\alpha } X(1))_+^{\alpha -1} ; X(1) \le 1 } 
\label{} \\
\le& \int_0^{\infty } \d t \e^{-qt} \{ |x|+t^{1/\alpha } \}^{\alpha -1} 
< \infty . 
\label{}
\end{align}
Since $ \bP_0(X(1) \ge y) \le cy^{-\alpha -1} $ for all $ y>1 $ for some constant $ c>0 $, 
we have 
\begin{align}
& \int_0^{\infty } \d t \e^{-qt} t^{1-\alpha } 
\bE_0 \sbra{ X(1)^{\alpha -1} - 1 ; X(1) > 1 } 
\label{} \\
=& \int_0^{\infty } \d t \e^{-qt} t^{1-\alpha } 
\int_1^{\infty } (\alpha -1) y^{\alpha -2} \bP_0(X(1) \ge y) \d y 
\label{} \\
\le& c(\alpha -1) 
\int_0^{\infty } \d t \e^{-qt} t^{1-\alpha } 
\int_1^{\infty } y^{-3} \d y < \infty . 
\label{}
\end{align}
Since $ (a+b)^p \le a^p+b^p $ for all $ a,b>0 $ and all $ 0<p<1 $, we obtain 
\begin{align}
& \int_0^{\infty } \d t \e^{-qt} 
\bE_0 \sbra{ \{ |x|+t^{1/\alpha } X(1) \}^{\alpha -1} ; X(1)>1 } 
\label{} \\
\le& \int_0^{\infty } \d t \e^{-qt} 
\bE_0 \sbra{ |x|^{\alpha -1} +t^{1-\alpha } X(1)^{\alpha -1} ; X(1)>1 } < \infty . 
\label{}
\end{align}
We thus obtain \eqref{eq: xvee0}. The proof is now complete. 
}

\subsection{Existence of duration density}

\begin{Lem} \label{lem: LA hol}
Suppose that {\bf (L1$ ' $)} and {\bf (L2)} are satisfied. 
Then the function $ z \mapsto r_z(0) $ 
can be analytically continuated to a holomorphic function on the domain 
$ \{ z \in \bC : \Re z>0 \} $, 
and the formula 
\begin{align}
r_z(0) = \frac{1}{\pi} \int_0^{\infty } 
\frac{z+\theta(\lambda)}{(z+\theta(\lambda))^2 + \omega(\lambda)^2} \d \lambda 
\label{eq: rz0}
\end{align}
is still valid on $ \{ z \in \bC : \Re z>0 \} $. 
\end{Lem}

\Proof{
For $ z \in \bC $ and $ \lambda>0 $, let 
\begin{align}
F_{\lambda}(z) = (z+\theta(\lambda))^2 + \omega(\lambda)^2 . 
\label{}
\end{align}
Let $ z=q+ix $ with $ q>0 $ and $ x \in \bR $. We have 
\begin{align}
& |F_{\lambda}(q+ix)|^2 
\label{} \\
=& \absol{ (q+ix+\theta(\lambda))^2 + \omega(\lambda)^2 }^2 
\label{} \\
=& \absol{ (q+\theta(\lambda))^2 + \omega(\lambda)^2 - x^2 + 2ix(q+\theta(\lambda)) }^2 
\label{} \\
=& \cbra{ (q+\theta(\lambda))^2 + \omega(\lambda)^2 - x^2 }^2 + \cbra{ 2x(q+\theta(\lambda)) }^2 
\label{} \\
=& x^4 + 2 \{ (q+\theta(\lambda))^2 - \omega(\lambda)^2 \} x^2 
+ \{ (q+\theta(\lambda))^2 + \omega(\lambda)^2 \}^2 
\label{} \\
=& 
(x^2 - \omega(\lambda)^2)^2 
+ (q+\theta(\lambda))^2 \{ (q+\theta(\lambda))^2 + 2 \omega(\lambda)^2 + 2 x^2 \} . 
\label{eq: LA hol1}
\end{align}

Note that 
\begin{align}
\partial _z \frac{z+\theta(\lambda)}{F_{\lambda}(z)} 
=& \frac{F_{\lambda}(z) - (z+\theta(\lambda)) \partial _z F_{\lambda}(z) }{F_{\lambda}(z)^2} 
\label{} \\
=& \frac{- (z+\theta(\lambda))^2 + \omega(\lambda)^2 }{F_{\lambda}(z)^2} . 
\label{eq: LA hol2}
\end{align}
If $ q>q_0>0 $, we obtain 
\begin{align}
\absol{ \partial _z \frac{z+\theta(\lambda)}{F_{\lambda}(z)} } 
\le& \frac{|z+\theta(\lambda)|^2 + \omega(\lambda)^2 }{|F_{\lambda}(z)|^2} 
\label{eq: LA hol3} \\
\le& \frac{(q+\theta(\lambda))^2 + x^2 + \omega(\lambda)^2 }{(q+\theta(\lambda))^2 \{ (q+\theta(\lambda))^2 + 2 \omega(\lambda)^2 + 2 x^2 \}} 
\label{} \\
\le& \frac{1}{(q+\theta(\lambda))^2 } 
\label{} \\
\le& \frac{1}{q_0 (q+\theta(\lambda)) } 
\in L^1(\d \lambda) . 
\label{}
\end{align}
Hence, by {\bf (L1$ ' $)}, we see that 
the right hand side of \eqref{eq: rz0} 
is holomorphic on $ \{ z \in \bC : \Re z>0 \} $. 
The proof is now complete. 
}

Let us give a sufficient condition for {\bf (T)} to be satisfied. 

\begin{Thm} \label{thm: LA rho}
Suppose there exist constants $ 1<\alpha <2 $, 
$ \un{c}_{\theta}>0 $ and $ \bar{c}_{\theta}>0 $ such that 
\begin{align}
\un{c}_{\theta} \lambda^{\alpha } \le \theta(\lambda) \le \bar{c}_{\theta} \lambda^{\alpha } 
, \quad \lambda>0 
\label{eq: LA rho ass1}
\end{align}
and $ \un{c}_{\omega}>0 $ and $ \bar{c}_{\omega}>0 $ such that 
\begin{align}
\un{c}_{\omega} \lambda^{\alpha } \le \omega(\lambda) \le \bar{c}_{\omega} \lambda^{\alpha } 
, \quad \lambda>0 . 
\label{eq: LA rho ass2}
\end{align}
Then it follows that the conditions {\bf (L1$ ' $)}, {\bf (L2)} and {\bf (T)} are satisfied. 
\end{Thm}

\Proof{
The condition {\bf (L1$ ' $)} is obviously satisfied. 
If the condition {\bf (L2)} were not satisfied, 
i.e., $ v=0 $ and $ \int_{(-1,1)} |x| \nu(\d x) < \infty $, 
then it would hold that 
\begin{align}
\frac{\theta(\lambda)}{\lambda} 
\le& \int_{(-1,1)} \frac{1-\cos \lambda x}{\lambda x} x \nu(\d x) 
+ \frac{4}{\lambda} \nu((-1,1)^c) 
\label{} \\
\to& 0 
\quad \text{as $ \lambda \to \infty $}, 
\label{}
\end{align}
which would contradict the assumption. Thus we see that {\bf (L2)} is satisfied. 

By Theorem \ref{thm: ent law} and by Lemma \ref{lem: LA hol}, we see that 
\begin{align}
\begin{split}
\int_0^{\infty } \e^{-zt} t^2 \vn(T_0 \in \d t) 
=& - (\partial _z)^2 \frac{1}{r_z(0)} 
\\
=& \frac{\partial _z^2 r_z(0)}{r_z(0)^2} 
- \frac{2 (\partial _z r_z(0))^2}{r_z(0)^3} 
\end{split}
\label{eq: rz1}
\end{align}
holds for $ \Re z>0 $. 
Once we prove that 
\begin{align}
\varphi(x) := \sbra{ 
\frac{\partial _z^2 r_z(0)}{r_z(0)^2} 
- \frac{2 (\partial _z r_z(0))^2}{r_z(0)^3} 
}_{z=1+ix} 
\in L^1(\d x) , 
\label{}
\end{align}
we may invert the Fourier transform so that 
\begin{align}
\int_0^{\infty } \e^{-ixt} \e^{-t} t^2 \vn(T_0 \in \d t) = \varphi(x) 
\label{}
\end{align}
and obtain 
$ \vn(T_0 \in \d t) = \rho(t) \d t $, where 
\begin{align}
\rho(t) = \frac{1}{\e^{-t} t^2} \cdot 
\frac{1}{2 \pi} \int_{-\infty }^{\infty } \e^{itx} \varphi(x) \d x . 
\label{}
\end{align}

Let us obtain estimates involving variables $ x \in \bR $ and $ \lambda>0 $. 
On one hand, using the assumption, we can easily obtain 
\begin{align}
c_1(1 + |x| + \lambda^{\alpha })^2 
\le |F_{\lambda}(1+ix)| 
\le c_2(1 + |x| + \lambda^{\alpha })^2 
\label{}
\end{align}
for some constants $ c_1>0 $ and $ c_2>0 $. 
We can also obtain 
\begin{align}
|r_{1+ix}(0)| 
\ge& | \Re r_{1+ix}(0) | 
\label{} \\
=& \frac{1}{\pi} \int_0^{\infty } \frac{\theta(\lambda)}{|F_{\lambda}(1+ix)|} \d \lambda 
\label{} \\
\ge& c_3 (1+|x|)^{\frac{1}{\alpha }-1} 
\label{}
\end{align}
for some constant $ c_3>0 $. 
By \eqref{eq: LA hol3} and \eqref{eq: LA hol1}, we have 
\begin{align}
\absol{ \sbra{ \partial _z \frac{z+\theta(\lambda)}{F_{\lambda}(z)} }_{z=1+ix} } 
\le& \frac{|1+ix+\theta(\lambda)|^2 + \omega(\lambda)^2 }{|F_{\lambda}(1+ix)|^2} 
\label{} \\
\le& \frac{1}{c_2 (1 + |x| + \lambda^{\alpha })^2} . 
\label{}
\end{align}
We now obtain 
\begin{align}
\absol{ \sbra{ \partial _z r_z(0) }_{z=1+ix} } 
\le& \frac{1}{\pi} \int_0^{\infty } 
\absol{ \sbra{ \partial _z \frac{z+\theta(\lambda)}{F_{\lambda}(z)} }_{z=1+ix} } \d \lambda 
\label{} \\
\le& c_4 (1+|x|)^{\frac{1}{\alpha }-2} 
\label{}
\end{align}
for some constant $ c_4>0 $. Hence we obtain 
\begin{align}
\absol{ \sbra{ \frac{2 (\partial _z r_z(0))^2}{r_z(0)^3} }_{z=1+ix} } 
\le& \frac{2 c_4^2 (1+|x|)^{\frac{2}{\alpha }-4}}{c_3^3 (1+|x|)^{\frac{3}{\alpha }-3}} 
\label{} \\
\le& c_5 (1+|x|)^{-\frac{1}{\alpha }-1} 
\label{}
\end{align}
for some constant $ c_5>0 $. 
On the other hand, since we have 
\begin{align}
\absol{ \sbra{ \partial ^2_z \frac{z+\theta(\lambda)}{F_{\lambda}(z)} }_{z=1+ix} } 
\le& 
\absol{ \sbra{ \frac{- 4 (z+\theta(\lambda)) \omega(\lambda)^2 }{F_{\lambda}(z)^3} }_{z=1+ix} } 
\label{} \\
\le& c_6 (1+|x|+\lambda^{\alpha })^{-3} 
\label{}
\end{align}
for some constant $ c_6>0 $, we obtain 
\begin{align}
\absol{ \sbra{ \partial^2 _z r_z(0) }_{z=1+ix} } 
\le& \frac{1}{\pi} \int_0^{\infty } 
\absol{ \sbra{ \partial^2 _z \frac{z+\theta(\lambda)}{F_{\lambda}(z)} }_{z=1+ix} } \d \lambda 
\label{} \\
\le& c_7 (1+|x|)^{\frac{1}{\alpha }-3} 
\label{}
\end{align}
for some constant $ c_7>0 $. We thus obtain 
\begin{align}
\absol{ \sbra{ \frac{\partial^2 _z r_z(0)}{r_z(0)^2} }_{z=1+ix} } 
\le& \frac{c_7(1+|x|)^{\frac{1}{\alpha }-3}}{c_3^2(1+|x|)^{\frac{2}{\alpha }-2} } 
\label{} \\
\le& c_8 (1+|x|)^{-\frac{1}{\alpha }-1} 
\label{}
\end{align}
for some constant $ c_8>0 $. 

Therefore we conclude that 
\begin{align}
|\varphi(x)| \le (c_5+c_8) (1+|x|)^{-\frac{1}{\alpha }-1} , 
\label{}
\end{align}
which proves that $ \varphi(x) \in L^1(\d x) $. 
The proof is now complete. 
}

\def\cprime{$'$} \def\cprime{$'$}

\end{document}